\documentclass[12pt]{article}

\newtheorem{theorem}{Theorem}[section]

\newtheorem{proposition}[theorem]{Proposition}

\newtheorem{example}[theorem]{Example}
\newtheorem{remark}[theorem]{Remark}

\newfam\msbfam
\font\tenmsb=msbm10  scaled \magstep1 \textfont\msbfam=\tenmsb
\font\sevenmsb=msbm7 scaled \magstep1 \scriptfont\msbfam=\sevenmsb
\font\fivemsb=msbm5  scaled \magstep1 \scriptscriptfont\msbfam=\fivemsb
\def\Bbb{\fam\msbfam \tenmsb}

\def\CC{{\Bbb C}}

\def\ss{\subseteq}
\def\ra{\rightarrow}

\def\lh{{\cal L}H}

 \def\HollowBox #1#2{{\dimen0=#1 \advance\dimen0 by -#2       
       \dimen1=#1 \advance\dimen1 by #2                       
        \vrule height #1 depth #2 width #2                    
        \vrule height 0pt depth #2 width #1                   
        \llap{\vrule height #1 depth -\dimen0 width \dimen1}%
       \hskip -#2                                             
       \vrule height #1 depth #2 width #2}}                   
 \def\BoxOpTwo{\mathord{\HollowBox{6pt}{.4pt}}\;}             

\def\endpf{\hfill $\BoxOpTwo$ \smallskip \\ }

\newfam\msbbfam
\font\tenmsbb=msbm10  scaled \magstep1 \textfont\msbbfam=\tenmsbb
\font\sevenmsbb=msbm7  scaled \magstep1 \scriptfont\msbbfam=\sevenmsbb
\font\fivemsbb=msbm5    scaled \magstep1 \scriptscriptfont\msbbfam=\fivemsbb

\usepackage{graphicx}

\usepackage{amsmath}

\usepackage{amssymb}

\usepackage{amsfonts}

\begin{document}

\begin{center}
\Large \bf Uniqueness Properties of Hardy Space Functions\footnote{{\bf Subject 
Classification Numbers:} 32A35, 32A50, 32A70, 30H10.}
\end{center}
\vspace*{.12in}
								   
\begin{center}
Steven G. Krantz
\end{center}
\vspace*{.2in}

\begin{center}
\today
\end{center}
\vspace*{.25in}

\begin{quote}
{\bf Abstract:}   We study boundary uniqueness properties of Hardy space functions in several
complex variables.  Along the way, we develop properties of the Lumer Hardy space.
\end{quote}
\vspace*{.25in}

\section{Introduction}

Ever since the early work of Riesz and Hardy, the Hardy spaces have
been an important cornerstone of the harmonic analysis of complex
variable theory.  Even so, the theory of Hardy spaces in several
complex variables is not nearly as well developed as one would like.

In particular, many of the basic results about boundary uniqueness
for these function spaces have not been rigorously and completely established (although
in some instances they are a part of the folklore).  It is
our purpose in this paper to set some of these matters straight.

Throughout this paper, if $\Omega \ss \mathbb{C}^n$ is a given bounded domain with $C^2$ boundary,
then we let $d\sigma$ denote the usual boundary area measure on $\partial \Omega$.
In other words, $d\sigma$ is $(2n-1)$-dimensional Hausdorff measure.
(See [KRA1, Appendix II].)  In situations when we are dealing with
a domain $\Omega_j$ or $\Omega_\epsilon$, we shall denote the boundary 
measure by $d\sigma_j$ or $d\sigma_\epsilon$.

A basic fact of life is that a Hardy space function $f$ on a
smoothly bounded domain $\Omega$ in $\mathbb{C}^n$ has
nontangential boundary limits (indeed {\it admissible}
boundary limits---see [STE1] and [KRA1, Ch.\ 8]) almost
everywhere with respect to the usual boundary area measure
$d\sigma$. Thus we may analyze the boundary trace function
$\widetilde{f}$ and pose questions about it. For example, if
the boundary trace function $\widetilde{f}$ vanishes on a set
of positive $\sigma$ measure, does it then follow that $f
\equiv 0$? Such a result is well known in one complex
variable, and follows from the canonical factorization (in
particular, from properties of the Nevanlinna class---see
below). But it is not generally established in several complex
variables.

Likewise, a biholomorphic mapping $\Phi: \Omega_1 \rightarrow \Omega_2$ of smoothly
bounded domains will have a boundary trace $\widetilde{\Phi}$.  
It follows either from the properness of the mapping (and its inverse)
or from the maximum principle that $\widetilde{\Phi}$ maps
$\partial \Omega_1$ to $\partial \Omega_2$.   And of course
$\widetilde{\Phi^{-1}}$ enjoys a similar property.

It is well known, in the one variable context, that
$\widetilde{\Phi}$ takes sets of zero boundary measure to sets
of zero boundary measure and sets of positive boundary measure
to sets of positive boundary measure (this follows, for
instance, because of Painlev\'{e}'s theorem\footnote{See [BUR]
for the provenance of the Painelv\'{e} result.} that the
biholomorphic mapping continues smoothly to the boundary, and
so does its inverse). Is an analogous result true in several
complex variables (in which context we do {\it not} have a
general boundary smoothness theorem)?

Interestingly, some of these questions relate to properties of Lumer's Hardy
spaces.  In the course of answering the above-raised questions, we also develop
some relevant properties of Lumer's spaces.  We also use Lumer's Hardy
spaces in some of the proofs.

It may be noted that the answers to some of these queries are straightforward if
it is known that biholomorphic mappings continue smoothly and univalently to the
boundary.  So, for example, we may affirmatively answer the second question for strongly
pseudoconvex domains and finite type domains (see [KRA1]).  But our purpose here
is to come up with statements and proofs that are valid for {\it all smoothly
bounded domains in $\mathbb{C}^n$}.  

It is a pleasure to thank E. Bedford, S. R. Bell, R. Burckel,
and J. A. Cima for helpful remarks regarding the work in this
paper.

\section{Basic Concepts and Notation}
      
In all of our discussions, a {\it domain} in $\mathbb{C}^n$ will be a connected, open set.
We will let $\Omega$ denote a bounded domain in $\mathbb{C}^n$ with $C^2$ boundary.  So
$$
\Omega = \{z \in \mathbb{C}^n: \rho(z) < 0\} \, ,
$$
with $\rho$ a given $C^2$ defining function satisfying $\nabla \rho \ne 0$ on $\partial \Omega$.

Recall, for $0 < p < \infty$, that
$$
H^p(\Omega) = \left \{f \ \hbox{holomorphic on $\Omega$}: \sup_{\epsilon > 0} \int_{\rho(\zeta) = - \epsilon} |f(\zeta)|^p \, d\sigma_\epsilon(\zeta)^{1/p} \equiv \|f\|_{H^p(\Omega)} < \infty \right \} \, ,
$$
where $d\sigma_\epsilon$ is boundary area measure on the boundary of $\Omega_\epsilon \equiv \{z: \rho(z) < - \epsilon\}$.
It is known (see [KRA1, Ch.\ 8]) that this definition is independent of the choice of defining
function $\rho$.  We also define $H^\infty(\Omega)$ to be the bounded holomorphic functions
on $\Omega$ with the obvious supremum norm.

We say that a function $g$ on a domain $\Omega$ has a {\it harmonic majorant} (resp.\ {\it pluriharmonic majorant}) $u$ 
if the harmonic (resp.\ pluriharmonic) function $u$ satisfies $|g(z)| \leq u(z)$ for all $z \in \Omega$.
It is a basic fact that a holomorphic function $f$ on $\Omega$
lies in $H^p(\Omega)$, $0 < p < \infty$, if and only if $|f|^p$ has a harmonic
majorant (see [STE1], [KRA1, Ch.\ 8]).

The trouble with $H^p$ as we have defined and discussed it here
is that it is not obviously a biholomorphically invariant
notion (as is, for instance, the concept of Bergman
space---see [KRA1, Ch. 1]). For this reason it is useful to
have the concept of Lumer Hardy space. We say that a
holomorphic function $f$ on $\Omega$ is in the {\it Lumer
Hardy $p$-space} ${\cal L}H^p$ if the function $|f|^p$ has a
{\it pluriharmonic majorant}. Since pluriharmonic functions
are preserved by biholomorphic mappings of several complex
variables, it is clear that the concept of Lumer Hardy space
is a biholomorphic invariant.

The norm on ${\cal L}H^p$ is specified as follows.
Fix a point $P_0 \in \Omega$.  Then $\|f\|_{{\cal L}H^p} \equiv \inf u^{1/p}(P_0)$, where
the infimum is taken over all possible pluriharmonic majorants $u$ of $|f|^p$.   
Note that different choices of $P_0$ give equivalent norms.

In spite of its favorable property of invariance, the Lumer space has a number of pathological
properties.  For instance (see [RUD]), the Lumer space $\lh^2$ is a Banach space but not a Hilbert
space.  

It is clear, just because a pluriharmonic function is obviously harmonic, that ${\cal L}H^p$ is
a subspace of $H^p$.  

Indeed, we have

\begin{proposition} \sl
The space $\lh^2(\Omega)$ in $\CC^n$, $n \geq 2$, is a Banach subspace of 
$H^2(\Omega)$ but is not a Hilbert space.
\end{proposition}
{\bf Proof:}  For simplicity take $\Omega$ to be the unit ball $B$.
Rudin [RUD] shows that $\lh^2$ contains a closed subspace that
is isomorphic to $\ell^\infty$.  Since every closed subspace of a Hilbert space
is a Hilbert space, it then follows that $\lh^2$ is {\it not} a Hilbert space.

To see that $\lh^2$ is a Banach space, suppose that $g_j$ is a
Cauchy sequence in $\lh^2$. Then each $|g_j|^2$ has a least
pluriharmonic majorant $u_j$ and $\{u_j(P_0)\}$ is Cauchy. But
the Harnack inequalities, together with a connectedness
argument, then show that the $u_j$ have a Cauchy subsequence
$\{u_{j_k}\}$ that converges uniformly on compact sets. It
follows that the $g_{j_k}$ converge. So the space is complete,
and hence is a Banach space. 
\endpf

\begin{proposition} \sl
The space ${\cal L}H^p$ is a proper subspace of $H^p$.
\end{proposition}
{\bf Proof:}  See [LUM] for a version of this result.

We restrict attention to $\Omega = B$, the unit ball in $\CC^2$, and to $p = 2$.
Set
$$
f(z) = (1 - z_1)^{-1} \, .
$$
Then, because 
$$
|z_2|^2 < 1 - |z_1|^2 \leq 2|1 - z_1| \, ,
$$
it is easy to see that $f \in H^2(B)$.  But it is clear that $|f|^2$ does not have a harmonic majorant
on the complex line $\{\zeta \langle 1, 0, \dots, 0 \rangle: \zeta \in \mathbb{C}\}$ and hence
does not have a pluriharmonic majorant on $B$.
\endpf 

\begin{proposition} \sl
The dual of $\lh^2(B)$ does not equal the dual of $H^2(B)$.
\end{proposition}
{\bf Proof:}	Of course this follows just by functional analysis
from Proposition 2.2.  But it is worthwhile to see the result explicitly.
			  
We restrict attention to the unit ball $B$ in complex dimension
$n = 2$.  If $g$ is holomorphic on
the unit ball in $\CC^2$, then
$$
g(z) = \sum_j g_j(z) \, ,
$$
where each $g_j$ is a homogeneous polynomial of degree $j$.  
Rudin [RUD] proved the following:  If $P \in \partial B$ and $j \geq 0$ is fixed, then
$$
|g_j(P)| \leq \|g\|_{\lh^2} \, .  \eqno (2.3.1)
$$

With this result in hand, we let $P = (1,0)$ and consider the linear functionals on $\lh^2$
defined by
$$
\varphi_n(f) = \frac{\partial^n}{\partial z_1^n} f(P) \ \ , \quad n \in \{1, 2, \dots \} \, .
$$  
This is a bounded linear functional by (2.3.1).

We apply the functional $\varphi_n$ to the $H^2$ function given by 
$$
f^k(z_1, z_2) = \frac{(z_1 + 1)^k}{2^k} 
$$
for $k \in \{1, 2, \dots \}$.  Of course this is a holomorphic peaking function at the boundary
point $P = (1, 0) \in \partial B$.  

Obviously $\varphi_n(f^k)= k (k - 1) \cdots (k - (n-1)) \cdot 2^{k-n}/2^k = k(k-1) \cdots (k - (n-1))/2^n$.
On the other hand, the value of
the least harmonic majorant of $|(f^k)^2|$ (which is simply the solution
of the Dirichlet problem with $|(f^k)|^2 \bigr |_{\partial B}$ as boundary
data) at 0 is $\leq C \cdot (1 + \alpha)^{2k}/2^{2k}$ for some $0 < \alpha < \sqrt{2} - 1$.  Here we have simply used the mean
value property for harmonic functions (or else we can see this assertion
by way of harmonic measure).  So the corresponding $H^2$ norm is $\leq C \cdot (1 + \alpha)^k/2^k$. And,
if $C \cdot (1+ \alpha)^k/2^k$ bounds $k(k - 1) \cdots (k - (n-1)) /2^n$ uniformly in $n$, then
$$
\frac{k(k-1) \cdots (k - (n-1))}{2^n} \leq C \cdot \frac{(1 + \alpha)^k}{2^k} 
$$
hence
$$
k(k-1) \cdots (k - (n-1)) \leq C \cdot (1 + \alpha)^k \cdot 2^{n-k} \, .
$$
We take $n = k/2$ (assuming that $k$ is even and $[k/2]$ otherwise)
and obtain
$$
k(k-1) \cdots (k - ((k/2) - 1)) < C \, .
$$
This is clearly a contradiction if $k$ is large.  So $\varphi_n$ is not bounded on $H^2$.
\endpf

It appears to be rather difficult to give a complete intrinsic
description of the dual of $\lh^2$.

\section{Boundary Regularity}

A first natural question to ask is whether
a biholomorphic mapping of domains extends to
an invertible mapping of the boundaries.

In one complex variable this question is fairly easy.
For suppose that $\phi: \Omega_1 \ra \Omega_2$ is a conformal
mapping of domains in $\CC$, each domain having a Jordan curve
as boundary.  Then a classical theorem of Carath\'{e}odory [GRK] tells
us that $\phi$ continues to a bicontinuous mapping of the closures.
The univalence follows easily.

In several complex variables we know, for domains with $C^2$ boundary
let us say, that a biholomorphic mapping $\Phi: \Omega_1 \ra \Omega_2$ 
has an almost-everywhere defined boundary mapping $\widetilde{\Phi}$ and
we also know that the inverse mapping $\Phi^{-1}$ has an almost-everywhere
defined boundary mapping $\widetilde{\Phi^{-1}}$.  But we do {\it not} know that 
$\widetilde{\Phi^{-1}} \circ \widetilde{\Phi} = \hbox{\rm id}$ or that
$\widetilde{\Phi} \circ \widetilde{\Phi^{-1}} = \hbox{\rm id}$.  The next
result addresses this matter.   

Even in one variable it appears that this matter is not well documented.  So 
we treat that case first.

\begin{theorem}  \sl
Let $\Omega_1$, $\Omega_2$ be bounded domains in $\CC$ with $C^2$ boundary.
Let $\Phi: \Omega_1 \ra \Omega_2$ be a conformal mapping.  Then the boundary
mappings $\widetilde{\Phi}$ and $\widetilde{\Phi^{-1}}$ are inverse to each other.
\end{theorem}
{\bf Proof:}  I thank S. R. Bell for helpful remarks about this proof.

First it must be noted that $\Phi$ and $\Phi^{-1}$ each extend
to be $C^1$ on the closures (see, for instance, [BEK], [POM]).  So we may think of $\Phi: \overline{\Omega_1} \ra \overline{\Omega_2}$
as a diffeomorphism and $\Phi^{-1}: \overline{\Omega_2} \ra \overline{\Omega_1}$ as a diffeomorphism.
It follows from the maximum principle, or from the properness of the mappings, that
$\Phi$ takes $\partial \Omega_1$ to $\partial \Omega_2$ and $\Phi^{-1}$ takes
$\partial \Omega_2$ to $\partial \Omega_1$.
		 
It follows from the Hopf lemma (see [KRA2]) that $\nabla \Phi$
is nonvanishing on $\partial \Omega_1$. The nonvanishing
derivative together with the smooth extension shows that each
boundary curve of $\Omega_1$ gets mapped to a particular
boundary curve of $\Omega_2$. And we see that, as a domain
point $p$ traverses a boundary curve in the domain $\partial
\Omega_1$, the corresponding image point $\Phi(p)$ in
$\partial \Omega_2$ traverses a boundary curve a certain
number of times. But the argument principle tells us that that
number of times is one. That proves the result. \endpf

\noindent {\bf Alternative View of Theorem 3.1:}  In case the boundaries
of the domains in question are not smooth, a result is still possible.
Consider for instance the case when $\partial \Omega_1$, $\partial \Omega_2$ have
Lipschitz boundary.  [We make this particular geometric hypothesis so that each boundary
point has a well defined concept of nontangential convergence, and also of radial convergence.]
Since each boundary curve (or curves) is Jordan, Carath\'{e}odory's theorem implies that
the conformal mapping extends continuously to the boundary, and so does its inverse.

Assume that $P \in \partial \Omega_1$ is a boundary point at
which $\Phi$ has a radial limit. If $\nu_P$ is the unit
outward normal vector at $P$, then we may consider $\lim_{t
\ra 0^+} \Phi(P - t \nu_P)$. This limit will certainly exist
(this is what we think of as the ``radial limit''), and
$t \mapsto \Phi(P - t \nu_P)$, $0 < t < \epsilon_0$, will describe some
curve in $\Omega_2$. And the fact that $\Phi^{-1} \circ \Phi =
\hbox{\rm id}$ implies that $\lim_{t \ra 0^+} \Phi^{-1} \circ
\Phi(P - t \nu_P)$ exists. But the generalized Lindel\"{o}f
principle proved in [LEV] then implies that the {\it
nontangential limit} of $\Phi^{-1}$ at $\widetilde{\Phi}(P)$ exists. 
From this we may conclude that the limit that defines $\widetilde{\Phi^{-1}}$
coincides with the pointwise boundary limit of $\Phi$.  Therefore $\widetilde{\Phi^{-1} \circ
\Phi} = \widetilde{\Phi^{-1}} \circ \widetilde{\Phi} =
\widetilde{\hbox{\rm id}} = \hbox{id}$. And that is what we
wish to conclude. \endpf

We note for the record that neither version of Theorem 3.1 obtains
(at least not immediately) in the context of several complex variables.
For there certainly is no theorem, except in special cases, that says
that biholomorphic mappings extend continuously to the boundary (or smoothly
to the boundary).  And the Lindel\"{o}f principle in several variables is
different from that in one complex variable (see [CIK], [KRA3]).

Now we turn to the several variable situation. Let $\Omega_1$,
$\Omega_2$ be bounded domains in $\CC^n$ with $C^2$ boundary.
Let $\Phi: \Omega_1 \ra \Omega_2$ be a biholomorphic mapping.
Then there is a $\sigma_1$-almost everywhere defined boundary
mapping $\widetilde{\Phi}: \partial \Omega_1 \ra \partial
\Omega_2$ and there is a $\sigma_2$-almost everywhere defined
boundary mapping $\widetilde{\Phi^{-1}}: \partial \Omega_2 \ra
\partial \Omega_1$. The mapping $\Phi$ converges to
$\widetilde{\Phi}$ both nontangentially and admissibly.
Likewise for the mappings $\Phi^{-1}$ and
$\widetilde{\Phi^{-1}}$.

The assertions about the existence of $\widetilde{\Phi}$ and
$\widetilde{\Phi^{-1}}$, and about the convergence of $\Phi$ to $\widetilde{\Phi}$ 
and the convergence of $\Phi^{-1}$ to $\widetilde{\Phi^{-1}}$ follow from
standard results about the boundary behavior of $H^\infty$ functions.
See [STE1] and [KRA1, Ch.\ 8].  
		   
\begin{theorem} \sl
The mapping $\widetilde{\Phi}$ and the mapping $\widetilde{\Phi^{-1}}$ 
are each one-to-one (almost everywhere) and onto (almost everywhere).
\end{theorem}
{\bf First Proof of Theorem 3.2:}   We need to show that $\widetilde{\Phi}$ 
is both one-to-one and onto in a suitable measure-theoretic
sense.  Let us begin with the surjectivity.  Seeking a contradiction, we suppose
that the image of $\widetilde{\Phi}$ misses a set $H \ss \partial \Omega_2$ of positive
$\sigma_2$ measure.  Let $f$ be the characteristic function of $H$.  For $\psi \in \lh^2(\Omega_1)$, we consider the linear functional
$$
\lambda: \psi \mapsto \int_{\partial \Omega_1} (f \circ \widetilde{\Phi})(\zeta) \cdot \psi(\zeta) \, d\sigma_1(\zeta) \, .
$$
Then it is immediate that
$\lambda$ is the zero functional.  But $f$ certainly does {\it not} induce the zero functional
on $\lh^2(\Omega_2)$.  And a biholomorphic mapping of domains will of course
induce an isomorphism of $\lh^2(\Omega_1)^*$ with $\lh^2(\Omega_2)^*$.  So that is a contradiction.
Thus $\widetilde{\Phi}$ is onto.

Next we treat the univalence.  This point is tricky because we must determine how to
formulate this univalence in an almost-everywhere sense.   Suppose that $f_1$ and
$f_2$ are continuous functions with disjoint compact supports on $\partial \Omega_1$ which
each map to the same function $g$ on $\partial \Omega_2$ in the sense that
$g \circ \widetilde{\Phi} = f_1$ on the domain of $f_1$ and $g \circ \widetilde{\Phi} = f_2$ on
the domain of $f_2$.   Consider the functionals
$$
\lambda_j(\varphi) = \int_{\partial \Omega_1} \varphi \circ \widetilde{\Phi}(\zeta) \cdot f_j(\zeta) d\sigma_1(\zeta)
$$
on $\lh^2(\Omega_2)$, $j = 1, 2$.  Now $\lambda_j$ may be thought of as the push-forward of the
functional on $\lh^2(\Omega_1)$ induced by $f_j$.  On the other hand, the Hahn-Banach
theorem tells us that the functional $\lambda_j$ extends to a linear functional
on $L^2(\partial \Omega_2)$ and is therefore given by inner product with an $L^2(\partial \Omega_2)$ function
$g_j$.  And $g_j$ is supported on the image of the support of $f_j$ under $\widetilde{\Phi}$.  In fact
$g_j$ must be $g$.  

But this says that $\lambda_1 = \lambda_2$.  And that tells us that the induced map
from $\lh^2(\Omega_1)^*$ to $\lh^2(\Omega_2)^*$ is not univalent.  That is false.

This contradiction completes the proof of univalence.
\endpf
\vspace*{.15in}

\noindent {\bf Second Proof of Theorem 3.2:} Here we prove a
weak version of the univalence. Restrict attention to $\CC^2$.
Assume that $\Omega_1 \ss \CC^2$ and $\Omega_2 \ss \CC^2$. 
Suppose that $E \ss \partial \Omega_1$ has positive $\sigma_1$
measure and that $E$ is mapped by $\widetilde{\Phi}$ to a
constant $\alpha \in \partial \Omega_2$. We will derive
therefrom a contradiction.

Let $P \in E$ be a point of $\sigma_1$ density. Using Fubini's
theorem, find a complex line $\ell$ transversal to $\partial
\Omega_1$ so that $\ell \cap \partial \Omega_1$ near $P$ is a
one-dimensional, smooth real curve and $\ell \cap E$ is a
subset of that curve having positive linear measure. We may
conclude that $\widetilde{\Phi}$ maps $\ell \cap E$ to
$\alpha$ (see [NAR]). But then we can look at a component
$\widetilde{\Phi}_j$ of $\widetilde{\Phi}$, $j = 1, 2$; it clearly maps
$\ell \cap E$ to $\alpha_j$. So we have a bounded,
scalar-valued holomorphic function of one complex variable
that is constant on a piece of the boundary having positive
measure. So $\widetilde{\Phi}_j \bigr |_{\ell \cap \Omega_1}$
is constant. Since we can make this argument for uncountably
many distinct complex lines $\ell$, we conclude that
$\widetilde{\Phi}_j$ is constant for each $j = 1, 2$. Thus $\Phi$ is
constant. That is of course impossible. We conclude that
$\widetilde{\Phi}$ is set-theoretically one-to-one (in a weak sense).

Now we wish to show that $\widetilde{\Phi}$ is set-theoretically onto.  Seeking a contradiction, we suppose
that the image of $\widetilde{\Phi}$ misses a set $F \ss \partial \Omega_2$ of positive $\sigma_2$ measure.
We can then, as above, find a transversal complex line $\ell$ so that $\ell \cap \partial \Omega_2$ is a one-dimensional, smooth
real curve and $\ell \cap F$ is a subset of that curve having positive linear measure.  
 
But then we see that $\widetilde{\Phi}^{-1}$ maps $\ell \cap F$ to the empty set.
That clearly contradicts our argument in the first two paragraphs of this proof.
\endpf

It should be noted that we have proved that $\widetilde{\Phi}$ is both
univalent and surjective (in a certain measure-theoretic sense). Of course
the same result holds true for $\widetilde{\Phi^{-1}}$. But we have {\it
not} proved that $(\widetilde{\Phi})^{-1} = \widetilde{\Phi^{-1}}$ nor
that $\widetilde{\Phi^{-1}} \circ \widetilde{\Phi} = \hbox{\rm id}$.
	     
\section{The First Main Result}

The theorem that we treat in this section is as follows.
Let $\Omega_1$, $\Omega_2$ be bounded domains in $\mathbb{C}^n$ with
$C^2$ boundary.  Let $d\sigma_j$ be the usual area measure on $\partial \Omega_j$, $j = 1, 2$.
Let $\Phi: \Omega_1 \rightarrow \Omega_2$ be biholomorphic.
Then $\Phi$ has boundary trace $\widetilde{\Phi}$ and $\Phi^{-1}$ has boundary
trace $\widetilde{\Phi^{-1}}$.  

\begin{theorem} \sl
Let $E \ss \partial \Omega_1$ have $\sigma_1$ measure
0 and $F \ss \partial \Omega_1$ have positive $\sigma_1$ measure.
Then $\widetilde{\Phi}(E)$ has $\sigma_2$ measure 0 and $\widetilde{\Phi}(F)$ has
positive $\sigma_2$ measure.
\end{theorem}
{\bf First Proof of Theorem 4.1:}  The first of these statements is the contrapositive of the second.  So we
concentrate on proving the second.


Seeking a contradiction, we supposed that $F \ss \Omega_1$ is a set of positive $\sigma_1$ measure and that 
$G \equiv \widetilde{\Phi}(F)$ has zero $\sigma_2$ measure.  Define $X$ to be the collection
of those functions which are characteristic functions of measurable sets in $\partial \Omega_2$ that
are disjoint from $G$.  And let $Y$ be the linear space generated by $\lh^2 (\Omega_2) \cup X$.
[Here we identify each element of $\lh^2(\Omega_2)$ with its boundary function.]

Consider the linear functional on $\lh^2 (\Omega_2)$ defined by
$$
\phi(f) = \int_{\partial \Omega_1}  \chi_F(\zeta) \cdot (f \circ \Phi)(\zeta) \, d\sigma_1(\zeta) \, .  \eqno (4.1.1)
$$
We may extend this functional to $X$ and then to $Y$ by simply setting $\phi(g) = 0$ for any $g \in X$ and then
extending by linearity.   Notice that this extended functional is also specified by the integral
in (4.1.1).

Now we may apply the Hahn-Banach theorem to extend the functional $\phi$ 
to a new functional $\widetilde{\phi}$ on all of $L^2(\partial\Omega_2)$.  So
of course $\widetilde{\phi}$ is given by integration against an $L^2(\partial \Omega_2)$ function $p(\zeta)$.  And,
since $\phi$ (and hence $\widetilde{\phi}$) annihilates $X$, it must be that $p$ lives on $G$.   But we know that $G$ has $\sigma_2$ measure 0.
So we are forced to conclude that the extended functional $\widetilde{\phi}$  is the identically zero functional.
But the functional
$$
\psi(h) = \int_{\partial \Omega_1} \chi_F(\zeta) \cdot h(\zeta) \, d\sigma_1(\zeta) 
$$
is certainly {\it not} the zero functional (simply take $h$ to be a polynomial with
positive real part).  

So we have determined that the canonical mapping induced by $\widetilde{\Phi}$ 
from the dual space of $\lh^2(\Omega_1)$ to the
dual space of $\lh^2(\Omega_2)$ sends a nonzero functional to the zero functional.
And that is impossible.
\endpf
\vspace*{.15in}

\noindent {\bf Second Proof of Theorem 4.1}  For convenience let us work in $\CC^2$.  Let $E \ss \partial \Omega_1$ have positive $\sigma_1$ measure.
We claim that $\widetilde{\Phi}(E)$ cannot have $\sigma_2$ measure zero.  Seeking a contradiction,
suppose not.

Let $P \in E$ be a point of $\sigma_1$ density. Choose a
complex line $\ell$ which is transversal to the boundary and
so that $\ell \cap \partial \Omega_1$ is a smooth curve and
further so that $\ell \cap E$ is a subset of that curve having
positive linear measure. It is easy to see, using Fubini's
theorem simultaneously in the domain and the range, that this
$\ell$ and the corresponding curve may be chosen so that the
image of $\ell \cap E$ under $\widetilde{\Phi}$ has measure 0.
Then $\Phi \bigr |_{\ell \cap \Omega_1}$ must be identically 0
by the one-variable theory. Since we can choose uncountably
many distinct complex lines $\ell$ of this nature, we may
conclude that $\Phi$ must be indentically 0. That is a
contradiction.

As usual, the statement that sets of measure 0 cannot get mapped to sets of positive measure
is just the contrapositive of what we have just proved.
\endpf

It would be incorrect to infer that a set of zero $\sigma$ measure in the boundary of a domain
{\it cannot} be a set of determinacy.  For instance, if a Hardy space function on
the bidisc vanishes on a set of positive 2-dimensional measure in the distinguished boundary
of the polydisc, then that function is identically zero.  See Theorem 5.1 below.

\begin{example}  \rm
The result of the theorem is false if the mapping is not biholomorphic.
As a simple example, let us first look at the one-complex-variable 
situation.  Let $D \ss \CC$ be the unit disc.  Recall the notable theorem of Bagemihl and Seidel [BAS]:
\end{example}

\begin{theorem} \sl
Let $E \ss \partial D$ be an $F_\sigma$ of first category.  Let $\varphi$ be any continuous
function on $D$.  Then there is a holomorphic function $f$ on $D$ so that
$$
\lim_{r \rightarrow 1^-} \biggl (\varphi(r \xi) - f(r \xi) \biggr ) = 0
$$
for every $\xi \in E$.
\end{theorem}
	     
Now let $E \ss \ \partial D$ be an $F_\sigma$ of first category and full measure.   Let $F \ss E$
have boundary measure 1 (note that $\partial D$ has boundary measure $2\pi$).  
Certainly there exists a continuous, complex-valued function $\varphi$ on $D$ with
values in the disc, finite unimodular radial boundary limit at every point of $E$ (call the boundary function $\widetilde{\varphi}$) and so that
$\widetilde{\varphi}(F)$ has measure 0.   Now the theorem of Bagemihl and Seidel guarantees the
existence of a holomorphic function $f$ on $D$ which has the same radial limits as $\varphi$ at each
point of $E$.  In particular, the boundary function $\widetilde{f}$ maps $F$ to a set of boundary
measure 0.

In the several complex variable context, consider $D^2 = D \times D$ and define
${\cal F}(z_1, z_2) = f(z_1)$ with $f$ as in the preceding paragraph.   Then this
${\cal F}$ takes the set $F \times D$ in the boundary of $D^2$ to a set of
measure zero.

Note that the paper [HAS] contains a version of the Bagemihl-Seidel theorem in
several complex variables that is valid on convex domains.  In particular, it			      
is valid on the unit ball in $\CC^n$.   They state their result for holomorphic
{\it functions}, but it is true as well for holomorphic mappings.  Thus, as in
the first paragraph, one can generate an example directly on the unit ball
of $\CC^n$.	

It should be noted that, in general, the functions provided by the Bagemihl-Seidel
theorem in [BAS] or in [HAS] are {\it not} bounded, and not necessarily in any $H^p$ class (nor in
the Nevanlinna class).  So, strictly speaking, the examples described in the three
paragraphs following Theorem 4.3 are not counterexamples to Theorem 4.1 in the 
case that $\Phi$ is not biholomorphic.
     
\section{Boundary Uniqueness}

First we recall the result in one complex variable.

\begin{theorem} \sl
Let $f$ be an $H^p$ function on the unit disc $D$, $0 < p < \infty$.
Let $\widetilde{f}$ be the boundary trace of $f$.  If $\widetilde{f}(\zeta) = 0$
for $\zeta$ in a set $E$ of positive linear measure in $\partial D$, then $f \equiv 0$.
\end{theorem}
{\bf Proof:}   Seeking a contradiction, we suppose that $f$ is {\it not}
identically equal to 0.

Certainly $f$ is in the Nevanlinna class $N$.  Therefore (see [DUR]),
$g \equiv 1/f$ also lies in $N$.  And obviously $\widetilde{g} = \infty$ on $E$.
But, setting $g_r(e^{i\theta}) = g(r e^{i\theta})$, we know that
$\log^+ |g_r| \ra \log^+ |\widetilde{g}|$ as $r \rightarrow 1^-$,
and that $\log^+ |\widetilde{g}|$ is
integrable.  On the other hand, $\log^+ |g_r| \ra \infty$ on the set $E$.  That
is a  contradiction.
\endpf

\noindent See [LEM, pp.\ 3, 17] for related results.

\begin{theorem} \sl
Let $f$ be a Hardy space function on the smoothly bounded 
domain $\Omega \ss \mathbb{C}^n$, and
suppose that the boundary trace function $\widetilde{f}$ vanishes on
a set $E \ss \partial \Omega$ of positive $\sigma$ measure.  Then $f \equiv 0$ on $\Omega$.
\end{theorem}
{\bf Proof:}  For simplicity of notation we take the complex dimension of the ambient
space to be 2.

By a simple change of coordinates, we may suppose that $P \in \partial \Omega$
is a point of $\sigma$-measure density of $E$ and that the unit outward normal direction to the boundary
at $P$ is ${\bf 1} = \langle 1, 0 \rangle$.  We may also suppose that $P$ is the origin
of coordinates.

Now let ${\bf d}$ be the complex tangential disc to $\partial \Omega$ at $P$ given by
$$
\zeta \longmapsto (1, \zeta) \equiv Z
$$
for $\zeta$ small (say $|\zeta| < \eta$).  For each such point $(1, \zeta) = Z$, we set
$$
D_Z = \Omega \cap \{Z + \xi {\bf 1}: \xi \in \mathbb{C}, |\xi| < \eta \} \, .
$$
Then each $D_Z$ is a complex disc that is transverse to the boundary $\partial \Omega$
and the collection of all $D_Z$ foliates a boundary neighborhood of $P \in \partial \Omega$.
Let $e_Z = D_Z \cap \partial \Omega$ for each $Z$.

We know that $f$ has (let us say) nontangential boundary limit
0 at each point of $E$. Hence, by Fubini's theorem, there are
an uncountable collection of $Z$ so that $f \bigr
|_{D_Z}$ has nontangential boundary limit 0 on a subset of $e_Z$
that has positive linear measure.   

Passing to a possibly smaller uncountable set of $Z$, we
can again apply Fubini's theorem to say that $f$ restricted to
each $D_Z$ is in $H^p$ of that $D_Z$. More precisely,
let $\epsilon_1 > \epsilon_2 > \cdots \ra 0$. Using our
earlier notation from Section 2, and invoking the fact that $f \in
H^p(\Omega)$, we know that $|f|^p$ is integrable on each
$\partial \Omega_{\epsilon_j} \cap D_Z$. Now using
standard results about equivalence of different $H^p$ norms
(see [KRA1, Ch.\ 8]), we may conclude that $|f|^p$ is integrable
on level sets of $\partial D_Z$, each $Z$. So $f \in
H^p(D_Z)$, each $Z$. Hence we may apply Theorem 5.1 on
each $D_Z$ to conclude that $f \bigr |_{D_Z} = 0$,
each $Z$.

But the structure of the zero set of a holomorphic function on a domain in $\mathbb{C}^2$ is well
known (see [KRA1, Ch.\  7].  It is not possible for a nonzero holomorphic function to vanish on
uncountably many disjoint leaves.  We can only conclude that $f \equiv 0$.
\endpf

One of the more interesting results of classical function theory is the following.

\begin{proposition} \sl
Let $f$ be a holomorphic function on the unit disc $D$ which has nontangential
boundary limit 0 on a set $E \ss \partial D$ of positive measure.  Then
$f \equiv 0$.
\end{proposition}

\begin{remark} \rm
Note that, in this result, we are {\it not} assuming that $f \in H^p$ for some $p$.
In fact there is {\it no}  growth condition on $f$.  But we {\it are} looking at {\it nontangential limits}
instead of radial limits.  This is a result of Privalov, for which
see [COL, pp.\ 145--146].
\end{remark}

\noindent {\bf Proof of the Proposition:}  By a standard measure-theoretic argument (see [STE2]), we may take the nontangential
approach regions to all have the same aperture.  Let $U$ be 
the domain obtained by taking the union of all the nontangential
approach regions which terminate at points of $E$, together with a suitable
simply connected, relatively compact subset of $D$.  Then $U$ is a simply connected domain with Lipschitz
boundary.   So there is a conformal mapping $\varphi$ of $U$ to the unit disc $D$.
But then $f \circ \varphi^{-1}$ is a holomorphic function on $D$ that extends continuously
to $\partial D$ and vanishes on a boundary set of positive measure.  

But now we are dealing with a {\it bounded} holomorphic function.  So of course
it must be that $f \circ \varphi^{-1} \equiv 0$ and so $f \equiv 0$.
\endpf

What is remarkable and notable is that the proposition is false if ``nontangential
boundary limit'' is replaced by ``radial boundary limit.'' This follows immediately
from the classical theorem of Bagemihl and Seidel [BAS] cited above.

Now we have a result of a similar nature in several complex variables.  For comparison
purposes, it is worth mentioning the theorem of Hakim and Sibony [HAS].

\begin{theorem} \sl
Let $\Omega \ss \CC^n$ be a smoothly bounded domain.  Let $f$ be a holomorphic
function on $\Omega$, and let $E \ss \partial \Omega$ be a set of positive
$\sigma$ measure.  If $f$ has admissible boundary limit 0 at each point of $E$ then
$f \equiv 0$.
\end{theorem} 
{\bf Proof:}  Near a point of $\sigma$-density of $E$, we may find a parametrized family of transversal
complex lines $\ell$, with the parametrizing set being real-two-dimensional and so that
{\bf (i)}  each $\ell \cap \partial \Omega$ is a smooth curve and {\bf (ii)}  each $\ell \cap E$
is a subset of $\ell \cap \partial \Omega$ having positive linear measure.  For each
such $\ell$, $f \big |_{\ell \cap \Omega}$ is a holomorphic function of one complex
variable that has nontangential boundary limit 0 on a set of positive linear measure.
So $f \bigr |_{\ell \cap \Omega} \equiv 0$.  Since this conclusion is true on each $\ell$,
we conclude that $f \equiv 0$ on $\Omega$.
\endpf

\section{Concluding Remarks}

The purpose of this paper has been to supply proofs of results that come up 
frequently in function-theoretic arguments of several complex variables, but
which are not well recorded in the literature.  We hope that the arguments and results
presented here will prove to be useful in other contexts.  Certainly the content of this
paper nicely complements the ideas in [KRA1, Ch.\ 8] and [STE1].

\newpage
									
\section*{References}
\vspace*{.25in}

\begin{enumerate}

\item[{\bf [BAS]}] F. Bagemihl and W. Seidel, Some boundary
properties of analytic functions, {\em Math. Zeitschr.}
61(1954), 186-199.
				       
\item[{\bf[BEK]}]  S. R. Bell and S. G. Krantz, Smoothness to
the boundary of conformal mappings, {\it Rocky Mountain J.\
Math.} 17(1987), 23--40.

\item[{\bf [BUR]}] R. Burckel, {\it An Introduction to
Classical Complex Analysis}, Vol. 1, Pure and Applied
Mathematics, 82, Academic Press, Inc., New York-London, 1979.

\item[{\bf[CIK]}] J. A. Cima and S. G. Krantz, A Lindel\"{o}f
principle and normal functions in several complex variables,
{\em Duke Math. Jour.} 50(1983), 303-328.

\item[{\bf [COL]}]  E. F. Collingwood and A. J. Lohwater,
{\it The Theory of Cluster Sets}, Cambridge University Press,
Cambridge, 1966.

\item[{\bf [DUR]}]  P. Duren, {\it The Theory of $H^p$
Spaces}, Academic Press, New York, 1970.

\item[{\bf [FEF]}] C. Fefferman, The Bergman kernel and
biholomorphic mappings of pseudoconvex domains, {\em Invent.
Math.} 26(1974), 1-65.

\item[{\bf [GRK]}]  R. E. Greene and S. G. Krantz, 
{\it Function Theory of One Complex Variable}, 3rd ed.,
American Mathematical Society, Providence, RI, 2006.

\item[{\bf [HAS]}] M. Hakim and N. Sibony, Boundary properties
of holomorphic functions in the ball of $\CC^n$, {\it Math.\
Ann.} 276(1987), 549--555.

\item[{\bf [HIR]}]  M. Hirsch, {\it Differential Topology},
Springer-Verlag, New York, 1976.

\item[{\bf [KRA1]}]  S. G. Krantz, {\it Function Theory of
Several Complex Variables}, 2nd ed., American Mathematical
Society, Providence, RI, 2001.

\item[{\bf [KRA2]}]  S. G. Krantz, {\it Partial Differential
Equations and Complex Analysis}, CRC Press, Boca Raton, FL, 1992.

\item[{\bf [KRA3]}] S. G. Krantz, The Lindel\"{o}f principle in
several complex variables, {\it Journal of Math.\ Anal.\ and
Applications} 326(2007), 1190--1198.

\item[{\bf[LEV]}]  O. Lehto and K. Virtanen, Boundary behavior and
normal meromorphic functions, {\em Acta Math.} 97(1957),
47-65.

\item[{\bf [LEM]}] L. Lempert, Boundary behavoiur of
meromorphic functions of several variables, {\it Acta
Mathematica} 144(1980), 1–-25.

\item[{\bf [LUM]}]  G. Lumer, Espaces de Hardy en plusiers
variables complexes, {\it C. R. Acad.\ Sci.\ Paris}, S\'{e}r.\
A--B 273(1971), 151--154.

\item[{\bf [NAR]}] A. Nagel and W. Rudin, Local boundary
behavior of bounded holomorphic functions, {\it Canadian
Journal of Mathematics} 30(1978), 583-–592.

\item[{\bf [POM]}] C. Pommerenke, {\it Boundary Behaviour of
Conformal Maps}, Grundlehren der Mathematischen
Wissenschaften, 299, Springer-Verlag, Berlin, 1992.

\item[{\bf [RUD]}]  W. Rudin, Lumer's Hardy spaces, {\it Mich.\
Jour.\ Math.} 24(1977), 1--5.

\item[{\bf [STE1]}]  E. M. Stein, {\it Boundary Behavior of
Holomorphic Functions of Several Complex Variables}, Princeton
University Press, Princeton, NJ, 1972.

\item[{\bf [STE2]}]  E. M. Stein, {\it Singular Integrals and
Differentiability Properties of Functions}, Princeton University
Press, Princeton, NJ, 1970.

\end{enumerate}
\vspace*{.25in}

\begin{quote} 
Steven G. Krantz \\ 
Department of Mathematics \\ 
Washington University in St. Louis \\ 
St.\ Louis, Missouri 63130 \\ 
{\tt sk@math.wustl.edu} 
\end{quote}

\end{document}